\documentstyle[11pt]{article}

\newtheorem{thm}{Theorem}[section]
\newtheorem{lem}[thm]{Lemma}

\newtheorem{pro}[thm]{Proposition}

\newtheorem{defi}[thm]{Definition}

\newcommand{\gm }{\Gamma }
\newcommand{\lon }{\longrightarrow }
\newcommand{\be }{\begin{eqnarray*}}
\newcommand{\ee }{\end{eqnarray*}}

\input amssym.def
\input amssym

\newcommand{\complex}{{\Bbb C}}

\newcommand{\frakg}{{\frak g}}
\newcommand{\frakh}{{\frak h}}

\newcommand{\half}{\frac{1}{2}}

\newcommand{\F}{{\cal F}}

\newcommand{\cald}{{\cal D}}

\newcommand{\smalcirc}{\mbox{\tiny{$\circ $}}}

\def\description label#1{\hfil\bf[#1]\hfil}
\newcommand{\ot}{\mbox{$\otimes$}}

\newcommand{\otr}{\ot_{R}}

\newcommand{\hopf}{(H, R, \alpha , \beta , m , \Delta , \epsilon )}
\newcommand{\eendr}{End_{k}R}
\newcommand{\phia}{\phi_{\alpha}}
\newcommand{\phib}{\phi_{\beta}}
\newcommand{\uqg}{U_{q}\frakg}
\newcommand{\otrf}{\otimes_{R_{\F}}}

\newcommand{\ea}{\mbox{$E_{\alpha}$}}
\newcommand{\eb}{\mbox{$E_{-\alpha}$}}

\def\sdp{\mathbin{\hbox{$\mapstochar\kern-.3333em\times$}}}
\def\pds{\mathbin{\hbox{$\times\kern-.55em\mapstochar\,$}}}

\newcommand{\wed}{\mathbin{\lower1.5pt\hbox{$\scriptstyle{\wedge}$}}}

\let\Tilde=\widetilde

\def\chigh{{\raise1.5pt\hbox{$\chi$}}}
\let\phi=\varphi
\def\til0{\Tilde{0}}

\def\dminus{\raise2pt\hbox{\vrule height1pt width 2ex}\hskip3pt}

\def\pback#1{\mathbin{{{\lower1.2ex\hbox{$\times$}}\atop #1}}}

\def\vlra{\hbox{$\,-\!\!\!-\!\!\!-\!\!\!-\!\!\!-\!\!\!
-\!\!\!-\!\!\!-\!\!\!-\!\!\!-\!\!\!\longrightarrow\,$}}

\def\gpd{\,\lower1pt\hbox{$\longrightarrow$}\hskip-.24in\raise2pt
             \hbox{$\longrightarrow$}\,}

\def\lgpd{\,\lower1pt\hbox{$\vlra$}\hskip-1.02in\raise2pt\hbox{$\vlra$}\,}

\def\llgpd{\,\lower1pt\hbox{$\vvlra$}\hskip-1.3in\raise2pt\hbox{$\vvlra$}\,
}

\begin{document}

\title{{\bf Quantum groupoids associated to  universal dynamical R-matrices
}}
\author{ PING XU \thanks{ Research partially supported by NSF
        grants   DMS97-04391.}\\
 Department of Mathematics\\The  Pennsylvania State University \\
University Park, PA 16802, USA\\
        {\sf email: ping@math.psu.edu }}

\date{}

\maketitle
\begin{abstract} By using twist construction, we
obtain a quantum groupoid  $\cald\ot_{q}\uqg$ for any  simple Lie algebra
$\frakg$. The underlying Hopf algebroid structure encodes
all the information of the corresponding elliptic
quantum group-the  quasi-Hopf algebras as
obtained by Fronsdal, Arnaudon et. al.  and Jimbo
et. al..
\end{abstract}

{\bf Groupo\"{\i}des quantiques  associ\'es \`a des R-matrices dynamiques
universelles}

{\bf R\'esum\'e}
En utlisant une construction twist\'ee, nous obtenons, pour toute alg\`ebre
de Lie simple $\frakg$, un groupo\"{\i}de quantique $\cald\ot_{q}\uqg$. La
structure d'alg\'ebro\"{\i}de de Hopf sous-jacente 
encode toute l'information contenue dans le groupe quantique elliptique
correspondant---cf. les quasi-alg\`ebres de 
Hopf obtenues par Fronsdal, Arnaudon et. al.  et Jimbo et. al..

{\bf Version fran\c caise abr\'eg\'ee}

Il y a, depuis peu, un int\'er\^et grandissant pour l'\'equation  de 
 Yang-Baxter quantique universelle, aussi connue sous le nom d'\'equation
de Geravais-Neveu-Felder~: 
\begin{equation}
\label{eq:dybe0}
R_{12}(\lambda+\hbar h^{(3)})R_{13}(\lambda )R_{23}(\lambda+\hbar h^{(1)})
=R_{23}(\lambda )R_{13}(\lambda +\hbar h^{(2)})R_{12}(\lambda ).
\end{equation} 

Ici, $R(\lambda )$ est une fonction m\'eromorphe de $\eta^*$ vers $\uqg \ot
\uqg$, o\`u $\uqg$ est 
un groupe quantique quasi-triangulaire, et $\eta \subset \frakg$ est une
sous-alg\`ebre de Cartan. 
Cette \'equation appara\^{\i}t naturellement dans de nombreux contextes de
la physique math\'ematique tels que 
le th\'eorie de Liouville quantique, l'\'equation de
Knizhnik-Zamolodchikov-Bernard quantique ou encore 
le mod\`ele de Caloger-Moser quantique. Une approche de cette \'equation
est celle Babelon et. al. \cite{BBB} qui 
utilise la th\'eorie des  quasi-alg\`ebres de Hopf. Consid\`erons une
fonction m\'eromorphe $F: \eta^* \lon \uqg \ot \uqg$ 
telle que $F(\lambda)$ est inversible pour tout $\lambda $. Posons
$R(\lambda )=F_{21}(\lambda )^{-1}R F_{12}(\lambda ), $ 
o\`u $R \in \uqg \ot \uqg$ est la $R$-matrice universelle standard pour le
groupe quantique $\uqg$. 
On peut v\'erifier que  $R(\lambda )$ satisfait l'\'equation
(\ref{eq:dybe0}) quand $F(\lambda )$ est de poids nul, 
et satisfait la condition de cocycle d\'ecal\'e suivante~:
 \begin{equation}
\label{eq:shifted0}
(\Delta_0  \ot  id )F(\lambda )  F_{12} (\lambda +\hbar h^{(3)})
 =  (id \ot  \Delta_{0} )  F (\lambda )F_{23}(\lambda  ),
\end{equation}
o\`u $\Delta_0$ est le co-produit dans $\uqg$.
Si nous supposons de plus que
\begin{equation}
\label{eq:co0}
 (\epsilon_{0} \ot id) F(\lambda )  =  1; \ \
(id \ot \epsilon_{0}) F (\lambda ) =  1,
\end{equation} o\`u $\epsilon_{0}$ est l'application co-unit\'e, on peut 
former un groupe quantique elliptique qui est une famille de
quasi-alg\`ebres de Hopf $(\uqg , \Delta_{\lambda})$ 
param\'etris\'ee par $\lambda \in \eta^*$:
$\Delta_{\lambda}=F(\lambda )^{-1}\Delta F (\lambda )$.

Le but de cet article est de montrer que l'\'equation (\ref{eq:shifted0})
d\'erive  aussi naturellement
de l'\'equation ``twistorielle" d'un groupo\"{\i}de quantique. Ceci conduit
\`a interpr\'eter un groupe quantique 
elliptique comme un groupo\"{\i}de quantique. Tout d'abord nous
introduisons une construction g\'en\'erale twist\'ee 
pour les alg\'ebro\"{\i}des de Hopf. De mani\`ere analogue au cas des
alg\`ebres de Hopf, nous montrons le \\\\
{\bf Th\'eor\`eme A}\ \ {\em 
Soient $\hopf$, un alg\'ebro\"{\i}de de Hopf d'ancre $\tau$, et $\F\in
H\otr H $ un twisteur. Alors, 
$(H, R_{\F}, \alpha_{\F} , \beta_{\F}, m, \Delta_{\F}, \epsilon)$
est un alg\'ebro\"{\i}de de Hopf admettant $\tau$ pour ancre. La
cat\'egorie mono\"{\i}dale de $H$-modules correspondante 
est \'equivalente \`a celle correspondant \`a $\hopf$.}

Maintenant, consid\'erons $H=\cald\ot \uqg [[\hbar ]]$ ($q=e^{\hbar }$),
o\`u $\cald$ est 
l'alg\`ebre des op\'erateurs diff\'erentiels sur $\eta^*$. $H$ est de
mani\`ere naturelle un alg\'ebro\"{\i}de de Hopf. 
Nous lui appliquons la construction twist\'ee~:\\\\
{\bf Th\'eor\`eme B} \ \  {\em Soit $F(\lambda )\in \uqg \ot \uqg$ une
fonction m\'eromorphe de poids nul. Alors, 
$\F =F(\lambda ) \Theta \in H\otr H$ est un twisteur (c.\`a.d. satisfait
l'\'equation (\ref{eq:cocycle})
et (\ref{eq:co-unit})) ssi les \'equations (\ref{eq:shifted0}) et
(\ref{eq:co0}) tiennent. Ici, $\Theta$ est d\'efini par l'\'equation 
(\ref{eq:theta}).}

Nous notons $\cald\ot_{q}\uqg$, le groupo\"{\i}de quantique obtenu en
twistant $\cald\ot \uqg [[\hbar ]]$ via $\F$. 
Comme con\-s\'e\-quen\-ce imm\'ediate du th\'eor\`eme A, nous avons le\\\\
{\bf Th\'eor\`eme C} \ \ {\em En tant que cat\'egorie mono\"{\i}dale, la
cat\'egorie des $\cald\ot_{q}\uqg$-modules 
est \'equivalente \`a celle des $\cald \ot \uqg [[\hbar ]]$-modules, et est
donc une cat\'egorie mono\"{\i}dale 
tress\'ee.}

Nous observons que $\cald\ot_{q}\uqg$ est co-associatif en tant
qu'alg\'ebro\"{\i}de de Hopf, bien que 
$(\uqg, \Delta_{\lambda})$ ne le soit pas. D\`es lors, la construction de
$\cald\ot_{q}\uqg$ consiste, d'une 
certaine mani\`ere, \`a \'etablir la co-associativit\'e en \'elargissant
l'alg\`ebre  $\uqg$ par produit 
tensoriel avec la partie dynamique $\cald$. Le lien entre ce groupo\"{\i}de
quantique et les  quasi-alg\`ebres de Hopf 
est, dans un certain sens, analogue \`a celui qu'il y a entre un fibr\'e
vectoriel et ses fibres.

\section{Introduction}
Recently, there is an increasing interest in the so called 
quantum universal  dynamical Yang-Baxter equation, also
known as Geravais-Neveu-Felder equation:
\begin{equation}
\label{eq:dybe}
R_{12}(\lambda+\hbar h^{(3)})R_{13}(\lambda )R_{23}(\lambda+\hbar h^{(1)})
=R_{23}(\lambda )R_{13}(\lambda +\hbar h^{(2)})R_{12}(\lambda ).
\end{equation}
 Here $R(\lambda )$ is a meromorphic function
from $\eta^*$ to $\uqg \ot \uqg$ ($q=e^{\hbar}$),
  $\frakg$  is a complex semi-simple Lie algebra
with Cartan subalgebra $\eta$, and $\uqg$ is a quasi-triangular
quantum group.
This equation  arises     
naturally from  various contexts in mathematical physics,
including quantum Liouville theory, quantum  
Knizhnik-Zamolodchikov-Bernard equation, and quantum
Caloger-Moser model.
 One   approach to this equation
due to Babelon et. al. \cite{BBB} is to use Drinfeld's
theory of quais-Hopf algebras. Consider a meromorphic function
$F: \eta^* \lon \uqg \ot \uqg$ such that $F(\lambda)$ is
invertible for all $\lambda$.  Set
$R(\lambda )=F_{21}(\lambda )^{-1}R F_{12}(\lambda ), $
where $R \in \uqg \ot \uqg$ is the standard universal $R$-matrix for the
quantum group $\uqg$. One can check that
$R(\lambda )$ satisfies Equation (\ref{eq:dybe}) if
$F(\lambda )$ is of zero weight,
and satisfies the following shifted cocycle condition:
\begin{equation}
\label{eq:shifted}
(\Delta_0  \ot  id )F(\lambda )  F_{12} (\lambda +\hbar h^{(3)})
 =  (id \ot  \Delta_{0} )  F (\lambda )F_{23}(\lambda  ),
\end{equation}
where $\Delta_0$ is the coproduct of $\uqg$.
If moreover we assume  that 
\begin{equation}
\label{eq:co}
 (\epsilon_{0} \ot id) F(\lambda )  =  1; \ \ 
(id \ot \epsilon_{0}) F (\lambda ) =  1, 
\end{equation}
where $\epsilon_{0}$ is the counit map, 
one can form an elliptic  quantum group, which is a family of
 quasi-Hopf algebras $(\uqg , \Delta_{\lambda})$
 parameterized by $\lambda \in \eta^*$:
$\Delta_{\lambda}=F(\lambda )^{-1}\Delta F (\lambda )$.
For $\frakg =\frak{sl}_{2}(\complex )$, a solution to Equations
(\ref{eq:shifted}) and (\ref{eq:co})
 was obtained by Babelon  \cite{Babelon}
in 1991.
For general simple Lie algebras, solutions were  recently
 found  independently by Arnaudon et. al. \cite{Arnaudon}  and Jimbo 
et. al. \cite{Jimbo} based on the approach of Fronsdal \cite{F}.

The purpose of this Note is to show that Equation (\ref{eq:shifted})
also arises naturally from the  ``twistor" equation
of a quantum groupoid. This 
 leads to  another interpretation
of an elliptic quantum group, namely 
 as a  quantum groupoid. The notion of quantum groupoids, as
a quantization of Lie bialgebroids,
was introduced in \cite{Xu} as a general framework
unifying quantum groups and star-products on Poisson
manifolds.
Roughly speaking, our construction goes as follows.
Instead of $\uqg$, we start with    the algebra
$H=\cald\ot \uqg [[\hbar ]]$ ($q=e^{\hbar }$), where
$\cald$ is the algebra of holomorphic differential operators on $\eta^*$.
$H$ is no longer a Hopf algebra. Instead it is a Hopf 
algebroid  as a  tensor product of  Hopf algebroids
  $\cald [[\hbar ]]$ and $\uqg$.
The  shifted  cocycle condition for $F(\lambda )$  is then
shown to be equivalent to the  twistor equation 
for this Hopf algebroid.
Using this twistor, we obtain 
 a  quantum groupoid  $\cald\ot_{q}\uqg$.
We note that $\cald\ot_{q}\uqg$ is co-associative as  a Hopf
algebroid, although $\Delta_{\lambda}$ is
not co-associative. The   construction of $\cald\ot_{q}\uqg$
is in some sense 
to restore  the co-associativity  by enlarging the algebra
 $\uqg$ by tensoring the dynamical part $\cald$.
The relation between this quantum groupoid and
   the quasi-Hopf algebras $(\uqg , \Delta_{\lambda})$
is, in a certain sense,  analogous to that between
 a fiber bundle and
its fibers. 
 We expect that this quantum groupoid will be
useful in understanding elliptic quantum groups, especially
their  representation theory.
The physical meaning of this quantum groupoid,   however,  still
needs to be explored.

The first part of the Note is  devoted to the introduction
of twist construction for Hopf algebroids.
 Then we retain to this particular example of universal
dynamical $R$-matrices in the second part.
We note that a different Hopf algebroid approach to the quantum
 dynamical Yang-Baxter equation
  was  also obtained by Etingof and
Varchenko \cite{EV}.

\section{Twist construction}
Let $\hopf$ be a Hopf algebroid (see \cite{Xu} for detail). 
By $\eendr$, we denote
the algebra of linear endmorphisms 
of $R$ (over the ground ring $k$). 
It is clear that $\eendr$ is an $(R,R)$-bimodule,
where $R$ acts on it from both  left and right  by  multiplications.
An {\em anchor} map is an $(R,R)$-bimodule map $\tau : H\lon \eendr$
satisfying certain axioms, which, roughly speaking, means
that $\tau$ is   a Hopf algebroid morphism. More precisely,
for $x\in  H$ and $a\in R$, we denote by $x(a)$ the
element in $R$ defined by $x(a)=\tau (x) (a)$,
and by $\phia$ and $\phib$ the maps $(H\otr H)\otimes R\lon H$
defined, respectively, by
$\phia (x\otr y\otimes a)=x(a)\cdot y$, and
$\phib (x\otr  y\otimes a)=x\cdot y(a)$.
Here $x, y\in H$, $a\in R$,  and the dot $\cdot$ denotes the 
$(R, R)$-action on $H$.  Note that $\phia$ and
$\phib$ are well defined since $\tau$ is an
$(R, R)$-bimodule map.
Then we require that
(i). $\tau : H\lon \eendr$ is an algebra homomorphism preserving
the identities;
(ii). $x(1_{R})=\epsilon x, \ \ \forall x\in H$;
(iii). $\phia (\Delta x \otimes a)=x\alpha (a)$ and \ $\phib (\Delta x
 \otimes a)=x\beta (a), \ \ \forall x\in H, a\in R$.

For $\uqg$, since $R= k$ and $\eendr \cong k$,
one can simply take  $\tau$ as
the counit. On the other hand,
if $H$ is the Hopf algebroid  $\cald$ of differential
operators over a manifold $M$ (see Example 2.1 in \cite{Xu}),
the usual action of differential operators on
functions is an anchor map.
More generally, if $H=\cald   \otimes
\uqg$ considered as a Hopf algebroid  over
$R=C^{\infty}(M)$, then  the map:
$\tau (D\otimes u)(f)=(\epsilon_{0} u) D(f), \ \ \forall D\in \cald ,
\ u\in \uqg , f\in C^{\infty}(M)$ is an anchor.
Here $\epsilon_0$ is the counit  of the Hopf algebra $\uqg$.

Given a Hopf algebroid  $\hopf$  with anchor $\tau$, let $\F$ be
an   element  in $ H\otr H$  satisfying:
\begin{eqnarray}
&& (\Delta \otr  id_{\scriptscriptstyle H} )\F  \cdot \F^{12} 
~ = ~ (id_{\scriptscriptstyle H} \otr  \Delta) \F  \cdot\F^{23} \ \ 
 \mbox{ in } \ \  H \otr H \otr H; \ \mbox{ and }  \label{eq:cocycle}\\
&& (\epsilon \otr id_{\scriptscriptstyle H}) \F ~ = ~ 1_{H}; \ \ 
(id_{\scriptscriptstyle H} \otr \epsilon) \F ~ = ~ 1_{H},
\label{eq:co-unit}
\end{eqnarray}
 where $\F^{12}=\F\otimes 1\in (H\otr H)\ot H$,   $\F^{23}=1\otimes \F
\in H\ot (H\otr H)$, and in Equation (\ref{eq:co-unit})
we have used the identification:  $R\otr H\cong H\otr R\cong H$.
 Define $\alpha_{\F}, \beta_{\F}:
R\lon H$, respectively, by
$\alpha_{\F}(a)=\phia (\F \ot a), \ \ \ \beta_{\F}(a)=\phib (\F \ot a),
\ \forall a\in R$. 
For any $a, b\in R$, set
$a*_{\F}b= \tau (\alpha_{\F}(a))(b)$. 
More explicitly, if $\F=\sum_{i}x_{i}\otr y_{i}$ for $x_{i}, y_{i}\in H$,
then 
$\alpha_{\F}(a)=\sum_{i} x_{i}(a) \cdot y_{i}, \ \beta_{\F}(a)=
\sum_{i} x_{i} \cdot y_{i}(a)$,  and $a*_{\F}b=\sum_{i} x_{i}(a) y_{i}(b),
\ \forall a, b\in R$.
Using Equations (\ref{eq:cocycle}) and (\ref{eq:co-unit}),
one can prove the following:
\begin{pro}
\begin{enumerate}
\item $(R, *_{\F})$ is an associative algebra, and
$1_{R}*_{\F}a=a*_{\F}1_{R}=a, \forall a \in R$.

\item $\alpha_{\F} :R_{\F}\lon H$ is an algebra homomorphism,
and $\beta_{\F} :R_{\F}\lon H$ is anti-homomorphism.
Here   $R_{\F}$ stands for the algebra $(R, *_{\F})$.
\end{enumerate}
\end{pro}

Moreover, it is not difficult to check that Equation (\ref{eq:cocycle})
also implies that
$\F   (\beta_{\F} (a)\otimes 1-1\ot \alpha_{\F} (a))=0 \ \mbox{ in } H\otr
H,
 \ \ \forall a\in R$.
As an immediate consequence, we have
 
\begin{lem}
\label{cor:tensor}
Let $M_{1}$ and $M_{2}$ be any $H$-modules. Then  $\F^{\#} (m_{1}\otrf m_{2})
=\F \cdot (m_{1}\ot m_{2}),\   m_{1}\in M_{1}, 
\mbox{ and } \ m_{2}\in M_{2}$,  is a well defined linear map.
Here and in the sequel, by $H$-modules
we always mean left $H$-modules.
\end{lem}

We say that  $\F $ is {\em invertible} if 
$\F^{\#}$ is a vector space isomorphism
for any $H$-modules $M_{1}$ and $M_{2}$. In this case, 
 in particular we can take $M_{1}=M_{2}=H$  so that  we  obtain
an isomorphism 
\begin{equation}
\label{eq:f}
\F^{\#} : H\otrf H \lon H \otr H .
\end{equation} 

\begin{defi}
An element $\F \in H\otr H$ is called a twistor if
it is invertible and  satisfies  both  Equations  (\ref{eq:cocycle}) and
(\ref{eq:co-unit}).
\end{defi}

Now assume that $\F$ is a twistor.  
Define a new coproduct  $\Delta_{\F}: H \lon H\otrf H $ by
$\Delta_{\F}=\F^{-1}\Delta \F$, which 
 means that
$\Delta_{\F}(x)={\F^{\#}}^{-1}
(\Delta (x)\F )$, $\forall x \in H$.\\\\
{\bf Theorem A}\ \ {\em Assume that $\hopf$ is a Hopf algebroid with anchor $\tau$, and
$\F\in H\otr H $ a twistor. Then
$(H, R_{\F}, \alpha_{\F} , \beta_{\F}, m, \Delta_{\F}, \epsilon)$
is a Hopf algebroid, which still admits $\tau$
as an  anchor. Its  corresponding monoidal category of 
$H$-modules is equivalent to that of  $\hopf$.}\\

We say that $(H, R_{\F}, \alpha_{\F} , \beta_{\F}, m, \Delta_{\F},
\epsilon)$
is obtained from $\hopf$ by twisting via $\F$.\\\\
{\bf Example 2.1} Let  $H=\cald (M) [[\hbar ]]$ be equipped with the
standard Hopf algebroid structure  over the base
  $R=C^{\infty}(M)[[\hbar]]$. As indicated earlier in this
section, $H$ admits
a natural anchor map.
A twistor $\F$ of the form: 
$\F= 1\otr 1+\hbar B_{1}+\cdots \in H\otr H [[\hbar ]]$ is equivalent to a
$*$-product on $M$.
The corresponding twisted Hopf algebroid was described
explicitly in \cite{Xu}.\\\\
{\bf Remark} If $\F_{1}\in H\otr H$ is a  twistor
to the Hopf algebroid $H$,
 and $\F_{2}\in
H\ot_{R_{\F_{1}}}H$ a twistor for the twisted
Hopf algebroid $H_{\F_{1}}$, then
the Hopf algebroid obtained by twisting $H$  via $\F_{1}$ then
via $\F_{2}$ is equivalent to that obtained by twisting
via $\F_{1} \F_{2}$. Here $\F_{1} \F_{2}\in H\otr H$ 
is understood as $\F_{1}^{\#}( \F_{2} )$, where
$\F_{1}^{\#} : H\ot_{R_{\F_{1}}}H\lon H\otr H$ is the map as
defined  in  Equation (\ref{eq:f}).

\section{Universal dynamical  R-matrices}
As in the introduction, $\frakg$ is a semi-simple Lie algebra over
$\complex$ with Cartan subalgebra $\eta$,  and
 $\uqg$ is a quasi-triangular quantum group.
In this section, $C(\eta^{*})$  always
denotes  the space of meromorphic functions on $\eta^{*}$, and
 $\cald$ the algebra
of holomorphic differential operators on $\eta^{*}$.
Let us consider $H=\cald\ot \uqg [[\hbar ]]$, where for
simplicity we assume that $q=e^{\hbar}$.
Then $H$ becomes  a Hopf algebroid  over the base algebra
$R=C(\eta^{*})[[\hbar ]]$ in a standard way.  Its
coproduct and counit  are denoted respectively by
$\Delta$ and $\epsilon$.
We  fix a basis of $\eta$, say $\{h_{1},  \cdots ,h_{k}\}$.
Let $\{\xi_{1}, \cdots , \xi_{k}\}$ be its dual basis.
This  defines a coordinate system $(\lambda_{1}, \cdots ,\lambda_{k})$
on $\eta^{*}$.

Set 
\begin{equation}
\label{eq:theta}
\theta =\sum_{i=1}^{k}(\frac{\partial}{\partial \lambda_{i}}\ot h_{i})
\in H\ot H , \ \mbox{ and } \Theta =\exp{\hbar \theta} \in H\ot H.
\end{equation}
Note that $\theta $, and hence $\Theta$, is independent of
the choice of the basis of $\eta$.
The following fact  can be easily verified.

\begin{lem}
\label{lem:theta}
$\Theta $ satisfies Equations (\ref{eq:cocycle}) and
(\ref{eq:co-unit}).
\end{lem}


In other words, $\Theta$ is  a twistor
 for  the Hopf algebroid $H$. As we shall see below,
it is this $\Theta$ that links  a 
shifted cocycle $F(\lambda)$ and   a Hopf algebroid twistor of $H$. \\\\
{\bf Theorem B}\ \ {\em
Assume that $F\in C(\eta^{*}, \uqg \ot \uqg)$ is of zero weight:
$[F(\lambda ), 1\ot h+h\ot 1]=0, \ \forall \lambda \in \eta^{*}, h\in \eta$.
Then $\F=F(\lambda )\Theta\in H\otr H$ 
is a twistor (i.e. satisfies Equations (\ref{eq:cocycle}) and
(\ref{eq:co-unit}))
 iff Equations (\ref{eq:shifted}) and (\ref{eq:co})  hold.}\\

Now assume that $F(\lambda )$ is  a solution to Equations 
 (\ref{eq:shifted}) and (\ref{eq:co})
so that we can form a quantum groupoid by twisting
$\cald\ot \uqg [[\hbar ]]$ via $\F$. 
The resulting 
 quantum groupoid is denoted by $\cald\ot_{q}\uqg$ ($q=e^{\hbar }).$

As an immediate consequence of  Theorem A, we have\\\\
{\bf Theorem C}\ \  {\em 
As a monoidal category, the category of   $\cald\ot_{q}\uqg$-modules
is equivalent to  that of $\cald \ot \uqg [[\hbar ]]$-modules,
and therefore is a braided monoidal category.}\\

One can describe $\cald\ot_{q}\uqg$  more explicitly.

\begin{pro}
\begin{enumerate}
\item $f*_{\F}g=fg, \ \ \forall f, g\in C(\eta^{*})$. I.e., $R_{\F}$ is
isomorphic to $C(\eta^{*})[[\hbar ]]$ with pointwise multiplication;
\item  $\alpha_{\F}f= \exp{(\hbar
\sum_{i=1}^{k}h_{i}\frac{\partial}{\partial \lambda_{i}}})f=\sum_{1\leq
i_{1}, \cdots  ,i_{n}\leq k}
\frac{\hbar^{n}}{n!}
 \frac{\partial^{n}f}{\partial \lambda_{i_{1}}\cdots
\lambda_{i_{n}}}h_{i_{1}}\cdots h_{i_{n}}$, $\forall  f \in C(\eta^{*})$;
\item  $\beta_{\F}f=f$,  $\forall  f \in C(\eta^{*})$.
\end{enumerate}
\end{pro}

To describe the relation between  $\Delta_{\F}$ and the 
quasi-Hopf algebra coproducts $\Delta_{\lambda}$.
We need a  ``projection" map from $H\otrf H$ to  $C(\eta^{*},  \uqg\ot \uqg
)$.
This can be defined as follows.
Let $Ad_{\Theta}: H\ot H\lon H\ot H$ be the adjoint operator:
$Ad_{\Theta}w=\Theta w \Theta^{-1}$,  $\forall w\in H\ot H$.
There exists  an obvious  projection  from $H\ot H$ to $C(\eta^{*}, 
\uqg\ot \uqg )$, which is just taking the 0th-order component.
Now composing $Ad_{\Theta}$ with this projection, one
obtains a map from $H\ot H$ to $C(\eta^{*},  \uqg\ot \uqg )$.
It is not difficult to see that this map descends to a
map $H\otrf H\lon C(\eta^{*},  \uqg\ot \uqg )$,
which is  denoted by $T$.

\begin{pro}
Let $i: \uqg\lon H$ be the natural inclusion. Then
$\Delta_{\lambda}=T\smalcirc \Delta_{\F}\smalcirc i$.
\end{pro}
{\bf Remark}  We  may replace
$\theta$  in  Equation (\ref{eq:theta}) by
$\tilde{\theta} = \half \sum_{i=1}^{k} (\frac{\partial}{\partial \lambda_{i}}
\ot h_{i} -h_{i}\ot \frac{\partial}{\partial \lambda_{i}} )
\in H\ot H $ and  set $\tilde{\Theta} =\exp{\hbar \tilde{\theta}} \in H\ot
H$.
One can show that $\tilde{\F}=F(\lambda )\tilde{\Theta }$
satisfies Equation (\ref{eq:cocycle} ) is equivalent to
 the following twisted cocycle condition for $F(\lambda )$:
$$(\Delta_0  \ot  id )F(\lambda )  F_{12} (\lambda +\half \hbar h^{(3)})
 =  (id \ot  \Delta_{0} )  F (\lambda )F_{23}(\lambda  -\half \hbar h^{(1)}
).$$
Then $R(\lambda )=F_{21}(\lambda )^{-1}R F_{12}(\lambda )$ satisfies 
\begin{equation}
\label{eq:dybe1}
R_{12}(\lambda+\half \hbar h^{(3)})R_{13}(\lambda -\half \hbar h^{(2)} )
R_{23}(\lambda+\half \hbar h^{(1)})
=R_{23}(\lambda  -\half \hbar h^{(1)})R_{13}(\lambda +\half\hbar h^{(2)})
R_{12}(\lambda  - \half \hbar h^{(3)}).
\end{equation}
In fact, both  $\Theta$ and $\tilde{\Theta}$ are obtained
from the deformation quantization of the canonical
symplectic structure  $T^{*}\eta^{*}$ using
the normal order and the Weyl order
respectively, so they are equivalent. This indicates that
solutions to Equation (\ref{eq:dybe}) and Equation (\ref{eq:dybe1})
should be  equivalent as well.

Finally,  when $F(\lambda )$ is the shifted cocycle obtained  in 
 \cite{Arnaudon}   \cite{Jimbo}, 
the classical limit (see \cite{Xu})
of $\cald\ot_{q}\uqg$ is the coboundary
Lie bialgebroid $(A, \Lambda )$, where $A=T \frakh^{*}\times \frakg$,
 $ \Lambda = \sum_{i=1}^{k} (\frac{\partial}{\partial \lambda_{i}}
\wedge  h_{i})+r(\lambda )\in \gm (\wedge^{2}A)$, and 
 $r(\lambda )$ is the standard dynamical $r$-matrix on $\frakg$:
$r(\lambda) \, = 
\half  \sum _{\alpha \in \Delta} \,
\coth ({\half } \ll \alpha, \lambda \gg) \ea \wedge \eb$.
Here $\ll \, , \, \gg$  is the Killing form of $\frakg$, $\Delta$ is
the set of roots of $\frakg$ with respect to $\eta$, the
$\ea$ and $\eb$'s are root vectors, and
$\coth (x)  =  {e^x + e^{-x} \over e^x - e^{-x}}$
is the hyperbolic cotangent function. 
As an consequence, we conclude that
this  coboundary Lie bialgebroid  is quantizable in
the sense of \cite{Xu}.

{\bf Acknowledgments.}  Most of
the work was completed during the author's visit in
IHES. He wishes to thank IHES
  for its   hospitality while this project was being done. 
 He is also grateful to  Melanie Bertelson and
 Pierre Bieliavsky for their help in
preparation of  the French abbreviated version.



\begin{thebibliography}{99}
\bibitem{Arnaudon}
Arnaudon, D., Buffenoir, E., Ragoucy, E., and
Roche, Ph.,
Universal solutions of quantum dymanical Yang-Baxter
equation,
{\em Lett. Math. Phys.} {\bf 44} (1998), 201-214.

\bibitem{Babelon}
Babelon, O., 
Universal exchange algebra for Bloch waves and Liouville theory,
{\em Comm. Math. Phys.} {\bf 139} (1991), 619-643.
\bibitem{BBB}
Babelon, O., Bernard, D. and Billey, E.,
A quasi-Hopf algebra interpretation of quantum
3j and 6j symbols and difference equations,
{\em Phys. Lett. B} {\bf 375} (1996), 89-97.

\bibitem{EV}
Etingof, P. and  Varchenko, A.,
Exchange dynamical quantum groups, math.QA/9801135.

\bibitem{Felder}
Felder, G.,
Conformal field theory and integrable systems associated
to elliptic curves,
{\em Proc. ICM, Zurich}, (1994), {\bf 2}, 1247-1255.

\bibitem{F}
Fronsdal, C.,
Quasi-Hopf deformation of quantum groups,
{\em Lett. Math. Phys.} {\bf 40} (1997), 117-134.

\bibitem{Jimbo}
Jimbo, M., Konno, H., Odake, and Shiraishi, J.,
Quasi-Hopf twistors for elliptic quantum groups,
q-alg/9712029


\bibitem{Xu}
Xu, P,
Quantum groupoids and deformation quantization,
{\em C. R.  Acad. Sci. Paris} Serie I.
{\bf  326} (1998), 289-294.

\end{thebibliography}
\end{document}